# Algebra transformations of the fundamental groups corresponding to those of Heegaard diagrams by the band moves

*By* Shunji HORIGUCHI

**Abstract.** This paper gives the basic result of [1](1997), i.e., a handle sliding and a band move of Heegaard diagrams correspond to a replacement and a substitution in relations of the fundamental groups derived from Heegaard diagrams, respectively (Theorem 12). Corollary 13 is a new addition for the homotopy 3-sphere.

**1. Preliminaries.**

Everything in this paper, we will be considering the piecewise linear point of view. $\partial X$, $Int(X)$, $Cl(X)$ shows the boundary, interior, closure of a point set $X$, respectively. Hereafter, notation $M^3$ denotes a closed, connected orientable 3-manifold unless otherwise stated.

In this section we give definitions. We begin with a definition of a handlebody.

**Definition 1.** Let $\{D_1,\cdots,D_n\}$ be mutually disjointed 2-disks and $h_i = D_i \times [0,1]$ $(i=1,\cdots,n)$. *A handlebody H of genus n* is a 3-ball(cube) $B^3$ with $n$ handles $\{h_1,\cdots,h_n\}$ so that the result of attaching $h_i$ with homeomorphisms throws $2n$ disks $D_i \times 0, D_i \times 1$ onto $2n$ disjointed 2-disks on $\partial B^3$. $H$ is represented as $B^3 + \bigcup_{i=1}^{n} h_i$ where $B^3 \cap h_i = \partial B^3 \cap \partial h_i = \{D_i \times 0, D_i \times 1\}$. A handlebody $H$ of genus $n$ is also called as *a solid torus of genus n*.

We note that $\partial H$ is an orientable or nonorientable closed surface of Euler characteristic $2-2n$ according as $H$ is orientable or nonorientable.

**Definition 2.** Let $H$ be a genus $n$ handlebody and $\{D_i\}(i=1,\cdots,n)$, mutually disjointed properly embedded 2-disks in $H$. If the $Cl(H-\{D_1\cup\cdots\cup D_n\})$ becomes a 3-ball, then the collection $\{D_i\}(i=1,\cdots,n)$ is called *a complete system of meridian disks of H* and $\{\partial D_i\}(i=1,\cdots,n)$ *a complete system of meridian circles of* $\partial H$.

Note that $\{D_1,\cdots,D_n\}$ cuts $\partial H$ into a 2-sphere with $2n$ holes.

**Definition 3.** (1) Let $H$ be an orientable genus $n(\geq 2)$ handlebody with the same presentation as in Def. 1. Fig. 1 shows two handles $h_i$ and $h_j$ of $H$. By an ambient isotopy of $H$, keeping $D_i \times 0$ fixed, and sliding $D_i \times 1$ along the direction of the line in $\partial(B^3 + h_j)$, $h_i$ goes over the $h_j$ and turns back to the first place. This operation is called *a handle sliding of $h_i$ about $h_j$*.

(2) Let $\{D_i\}(i=1,\cdots,n)$ be a complete system of meridian disks of $H$ and $m_i(=\partial D_i)$ a complete system of meridian circles of $\partial H$. Let $\alpha$ be an arc on $\partial H$ that joins two chosen meridians $m_i, m_j$ and $Int(\alpha) \cap (m_i \cup m_j) = \phi$. See Fig. 2. Let $N(m_i + \alpha + m_j, \partial H)$ be a regular neighborhood of $m_i + \alpha + m_j$ on $\partial H$. $\partial N$ consist of three circles. Out of the three circles, two are isotopic to $m_i, m_j$ and then the remainder is not isotopic to them. Let the notation of remainder be $m_{ij}$. $m_{ij}$ is called *a band sum of $m_i$ and $m_j$ (with respect to the band $\alpha$)*. It has also the very pleasant property that bounds a disk and it is homeomorphic to $D_i$ and $D_j$. Changing the label $m_{ij}$ into $m_i$ ($m_j$ resp.) is called *a band move of $m_i$ ($m_j$ resp.)*.

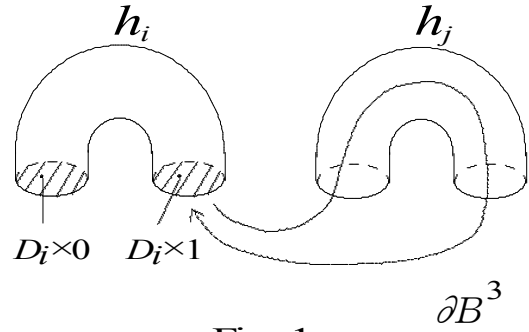

Fig. 1

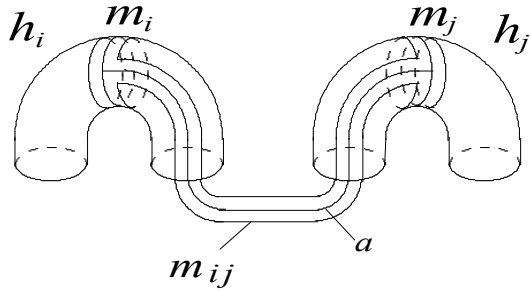

Fig. 2

**Definition 4.** A closed, connected 3-manifold $M^3$ is represented with a union of two handlebodies $H_1$, $H_2$ along their boundaries in $M^3$; $M^3 = H_1 \cup H_2$ so that $H_1 \cap H_2 = \partial H_1 \cap \partial H_2 = \partial H_1 = \partial H_2$. $\partial H_1 (= \partial H_2)$ is a closed surface of genus $n(\geqq 0)$. Let the surface be $F$. $H_1$($H_2$ resp.) and $F$ are orientable or nonorientable according as $M^3$ is orientable or nonorientable. A triplet $(H_1, H_2, F)$ or $M^3 = H_1 \cup H_2$ is called *a Heegaard splitting of $M^3$ with genus $n$* and $H_1$($H_2$ resp.), *a Heegaard-handlebody*. $F$ is called *a Heegaard-surface* and the integer $n(\geqq 0)$, *Heegaard genus*. Let $U$ and $V$ be disjointed handlebodies with the same genus. Let $f: U \to V$ be a homeomorphism so that $f|\partial U : \partial U \to \partial V$ is an orientation-reversing homeomorphism. Gluing together $\partial U$ of $U$ and $\partial V$ of $V$ by $f$, we obtain $M^3$. Then $M^3$ is denoted as $(M^3; U, V, f)$. It is called *a genus $n$ Heegaard splitting of $M^3$ concerning $f$*. In $(M^3; U, V, f)$, by replacing $f^{-1}(V)$ with $V$, one can regard $(M^3; U, V, f)$ as $(U, V, F)$ of $M^3$.

**Definition 5.** Suppose $(H_1, H_2, F)$ is a genus $n(\geqq 1)$ Heegaard splitting of $M^3$. Let $\{D_1, \cdots, D_n\}, \{D_1', \cdots, D_n'\}$ be a complete system of meridian disks of $H_1, H_2$, respectively. Let $\{m\} = \{m_1, \cdots, m_n\} = \{\partial D_1, \cdots, \partial D_n\}$, $\{l\} = \{l_1, \cdots, l_n\} = \{\partial D_1', \cdots, \partial D_n'\}$. Then $(H_1; m, l)$ $((H_2; l, m)$ resp.) is called *a genus $n$ Heegaard diagram associated with $(H_1, H_2, F)$*.
$(m, l)((l, m)$ resp.) are called *meridian-longitude* systems of $(H_1; m, l)((H_2; l, m)$ resp.).

**Definition 6.** By an ambient isotopy of $H$, a genus $n(\geqq 1)$ handlebody $H$ is deformed such as shown in Fig. 3.

If a genus $n$ Heegaard diagram $(H_1;m,l)$ satisfies the conditions of $m_i \cap l_j = \{a\ point\}\ (i = j)$ and $m_i \cap l_j = \phi\ (i \neq j)$, then $(H_1;m,l)$ is called *a canonical genus $n$ Heegaard diagram* of the 3-sphere.

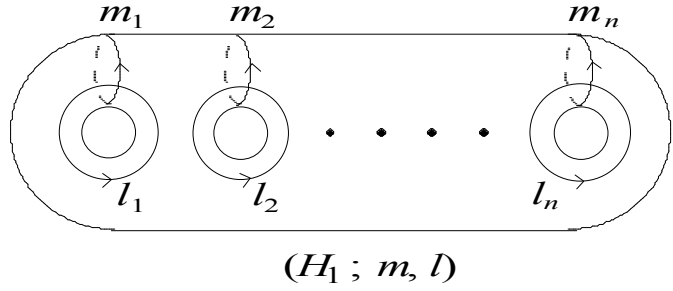

$(H_1 ; m, l)$

Fig. 3

Let $(H_1; m_1, \cdots, m_n, l_1, \cdots, l_n)$ be a genus $n$ Heegaard diagram associated with $(H_1, H_2, F)$ of $M^3$. We may assume that $(m_1 \cup \cdots \cup m_n) \cap (l_1 \cup \cdots \cup l_n)$ consists at most of finite points (by an argument of general position).

**Definition 7.** The number of finite points of $\{m\} \cap \{l\} = (m_1 \cup \cdots \cup m_n) \cap (l_1 \cup \cdots \cup l_n)$ is called *a number of cross points with* $(H_1;m,l)$ or $(H_2;l,m)$.

## 2. Transformations of Heegaard diagrams.

We begin with an obvious Proposition.

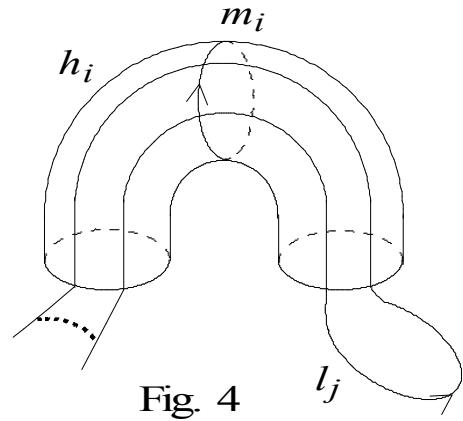

**Proposition 8.** *Let Fig. 4 be a part of Heegaard diagram $(U;m,l)$. The longitude $l_j$ crosses the meridian $m_i$, turns back to $m_i$ and crosses $m_i$ again. Then, there exists a transformation of $(U;m,l)$ so that a part of $l_j$ deforms to the dotted line and it does not cross $m_i$. It does not change the Heegaard genus but decreases the number of cross points, as many as* 2.

Fig. 4

**Definition 9.** The above transformation is called *a canceling for a Heegaard diagram*.

If the diagram like Fig. 4 appears, then we always do the above correction.

Let the following figure U-A be a part of Heegaard diagram $(U;m,l)$. The longitudes $\{l_{ij_1}, \cdots, l_{ij_l}\}\ (l \geqq 0)$ go around side by side on the two handles $h_i$ and $h_j$. The longitudes $\{l_{i_1}, \cdots, l_{i_p}\}, \{l_{j_1}, \cdots, l_{j_q}\}$ go around on $h_i, h_j$, respectively. It shows the general case that longitudes run on handles $h_i$ and $h_j$. In a special case that a character $l$ on the lower right equals to 0, there are not longitudes that run on $h_i$ and $h_j$. V-A is a part of $(V;l,m)$, and the dual part of U-A. The longitude $m_i, m_j$ crosses the meridians $\{l_{i_1}, \cdots, l_{i_p}, l_{ij_1}, \cdots, l_{ij_l}\}, \{l_{j_1}, \cdots, l_{j_q}, l_{ij_1}, \cdots, l_{ij_l}\}$ ,

respectively.

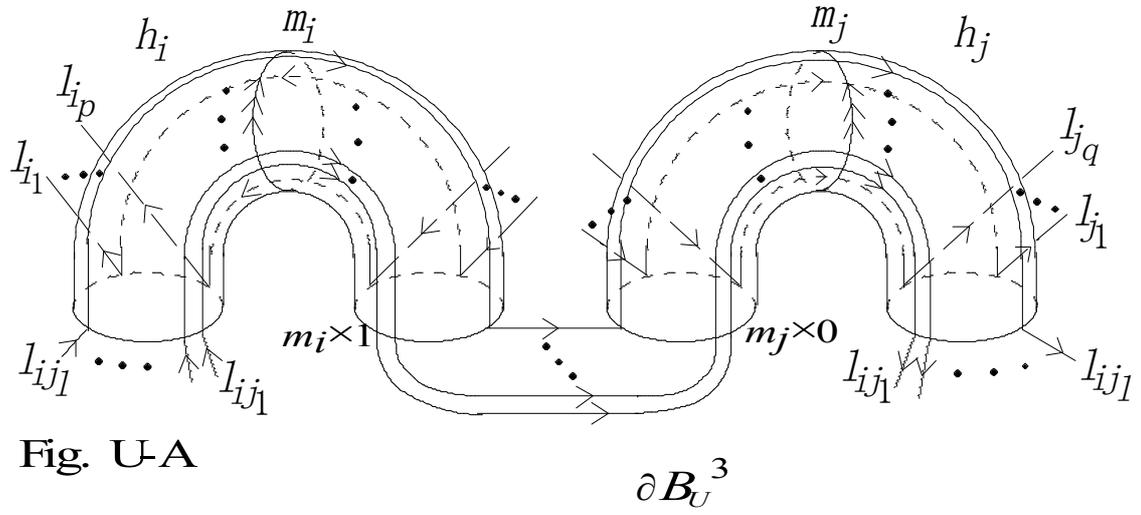

**Fig. U-A**

$\partial B_U{}^3$

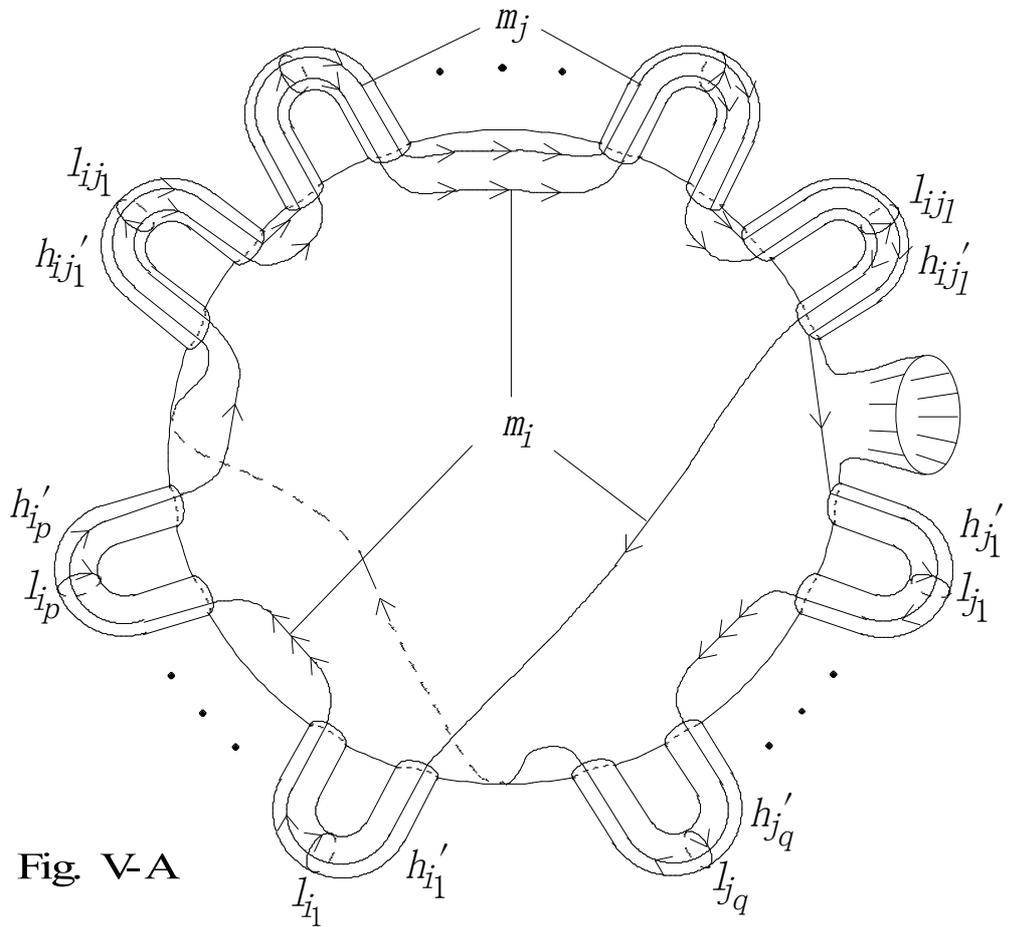

**Fig. V-A**

By the handle sliding of $h_i$ about $h_j$ along the directions of the longitudes $\{l_{ij_1},\cdots,l_{ij_l}\}$ in $\partial(B_U{}^3 + h_j)$, U-B is obtained from U-A.

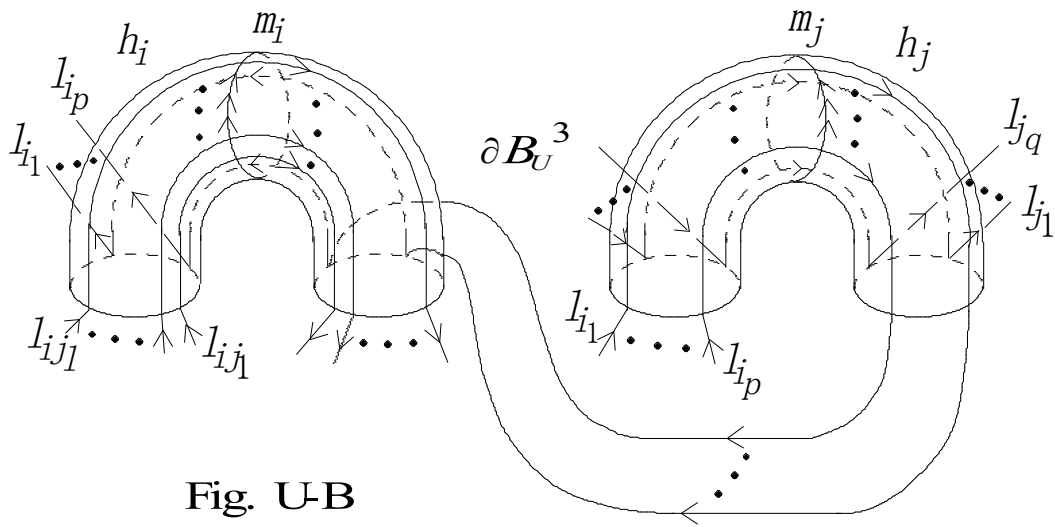

Fig. U-B

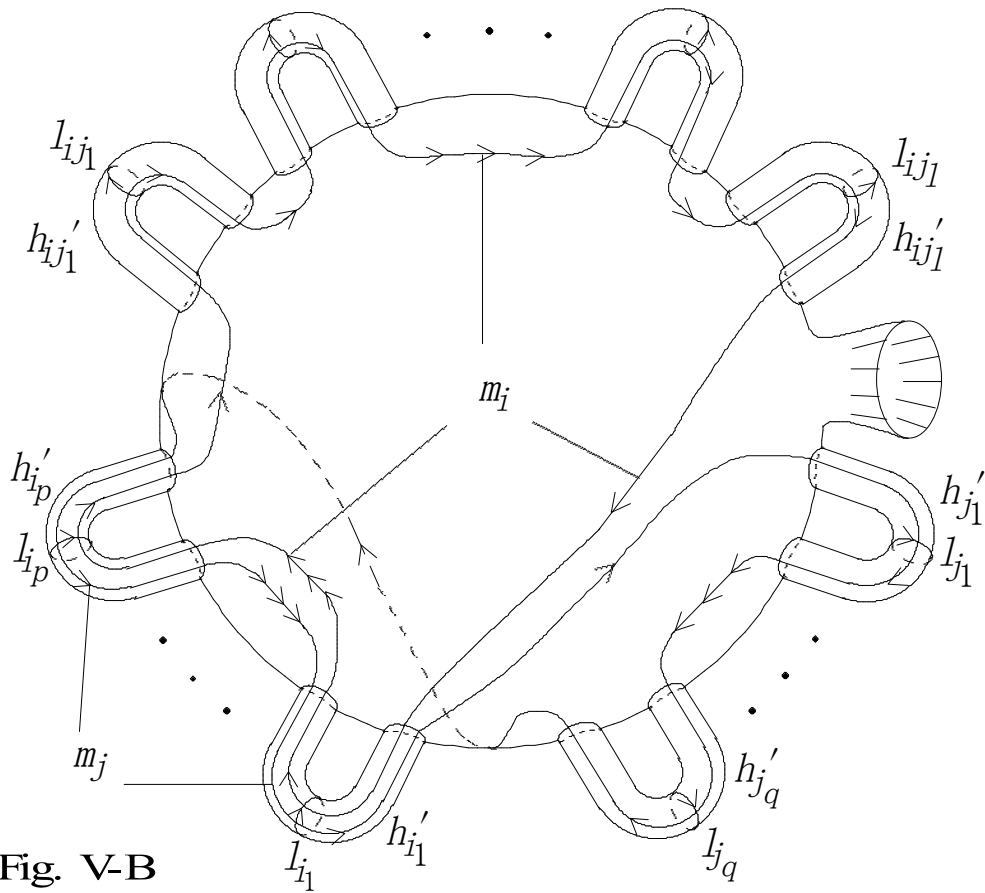

Fig. V-B

This handle sliding is the same as a band sum of $m_i$ and $m_j$ with respect to a part of a longitude $l_{ij_k}$ $(1 \leq k \leq l)$ that exists between $m_i$ and $m_j$, and a band move of $m_j$. In U-B, $\{l_{ij_1}, \cdots, l_{ij_l}\}$ go around on $h_i$ (not on $h_j$), $\{l_{i_1}, \cdots, l_{i_p}\}$ go around on both the $h_i$ and $h_j$, and $\{l_{j_1}, \cdots, l_{j_q}\}$ do not change the way of running. The dual transformation from V-A into V-B means that the band move in U-A is carry out in V-A. That is, each meridian $l_{ij_k}$ is cut into two segments by the two longitudes $m_i$ and $m_j$. Let the shorter segment be $\alpha$. From $m_i + \alpha + m_j$, we may construct a band sum $m_{ij}$ of $m_i$ and $m_j$, and carry out a band move of $m_j$. Next by an ambient isotopy and reorienting $m_j$, V-B is obtained. It is also obtained by handle slidings of $h_{j_1}'$ about $\{h_{ij_l}', h_{ij_{l-1}}', \cdots, h_{ij_1}', h_{i_p}', h_{i_{p-1}}', \cdots, h_{i_1}'\}$. In V-B, $m_i$ does not change the way of running, and here $m_j$ comes to cross $\{l_{j_1}, \cdots, l_{j_q}, l_{i_p}, \cdots, l_{i_1}\}$. The case of ($l=0$), if we can draw a band $\beta$ which reaches to $m_j \times 0$ via $m_i \times 1$ as it does not intersect the longitudes, then we can handle sliding $h_i$ about $h_j$ along $\beta + \partial h_j$. However, this only obtains more complex Heegaard diagram.

In like manners, a handle sliding of $h_j$ about $h_i$, and a band move of $m_i$ are obtained.

By cancelings for a Heegaard diagram in Def. 9, note that if $m \cap l \neq \phi$, then longitudes that run between handles $h_i$ and $h_j$ is only a type of Fig. U-A.

By the above transformation we have;

**Theorem 10.** *The transformation from U-A into U-B is carry out by a handle sliding of $h_i$ about $h_j$, or a band move of $m_j$. The dual transformation from V-A into V-B is carry out by a band move of $m_j$. These transformations do not change the Heegaard genus but change the number of cross points as many as $|l-p|$.*

In U-A(V-A resp.), we can carry out a band move for two meridians(two longitudes resp.) $m_i$, $m_j$ in Theorem 10.

### 3. Transformations of the fundamental groups.

To state our result precisely, we prepare algebra calculations for groups.

**Definition 11.** Let $\langle a_1, \cdots, a_n \mid r_1 = 1, \cdots, r_m = 1 \rangle$ denotes a presentation of a finitely generated group, where $a_1, \cdots, a_n$ are generators and relator $r_i$ is a word in the $a_i^\varepsilon$'s $(\varepsilon = \pm 1)$. We underline to the letters which are operated.

*Replacements letters*; if there are relations $\underline{a_i^\varepsilon a_j^\varepsilon} w_k = 1 (k=1,\cdots,\alpha)$ where $w_k$ is a word in the $a_i^\varepsilon$'s $(\varepsilon = \pm 1)$, then replace the generator $a_i$, letters $a_i^\varepsilon a_j^\varepsilon$ by a new letter $\tilde{a}_i$ (this becomes a new generator).

*Substitution*; if there are two relations $w_1 \underline{a_{i_1}^\varepsilon \cdots a_{i_\alpha}^\varepsilon} = 1$ and $w_2 \underline{a_{i_1}^\varepsilon \cdots a_{i_\alpha}^\varepsilon}$ where $a_{i_k} (k=1,\cdots,\alpha)$ is a generator and $a_{i_1}^\varepsilon \cdots a_{i_\alpha}^\varepsilon$ is a common word, then substitute $a_{i_1}^\varepsilon \cdots a_{i_\alpha}^\varepsilon = w_1^{-1}$ for $a_{i_1}^\varepsilon \cdots a_{i_\alpha}^\varepsilon$ of $w_2 a_{i_1}^\varepsilon \cdots a_{i_\alpha}^\varepsilon$.

Each above algebra calculation preserves isomorphism of a group.

Let $(U,V,F)$ be a genus $n(\geq 1)$ Heegaard splitting of $M^3$ and $(U;m,l)$ a Heegaard diagram of $(U,V,F)$. $\{m\} = \{m_1,\cdots,m_n\}$ and $\{l_1,\cdots,l_n\}$ are meridian-longitude systems. Let each $m_i, l_i$ be oriented. By applying Van Kampen's theorem to $U \cup V$, we may obtain a well-known presentation of a fundamental group $\pi_1(M^3)$;

$$\pi_1(M^3) = \langle m_1,\cdots,m_n \mid \hat{l}_1 = 1,\cdots,\hat{l}_n = 1 \rangle \quad (1)$$

We read that $m_1,\cdots,m_n$ are regarded as the generators of the meridians $m_1,\cdots,m_n$ and the relator $\hat{l}_j$ is a word in the $m_i^{\pm 1}$'s obtained by running once around the $l_j$, i.e., while we take a turn round $l_j$ according to the orientation of $l_j$, we read the label $m_i$ continuously as $m_i^{+1}$ ($m_i^{-1}$ resp.) if $l_j$ crosses $m_i$ from the left side (the right side resp.) to the right side (the left side resp.) of $m_i$. See Fig. 5. In the relator $\hat{l}_j$, we may start reading from any $m_i$ in $\hat{l}_j$ because the word $\hat{l}_j$ becomes a cyclic word by joining both ends of $\hat{l}_j$ and preserving the sequential order of letters in $\hat{l}_j$. Therefore $\hat{l}_j$ is uniquely defined up to cyclic permutations and inversions.

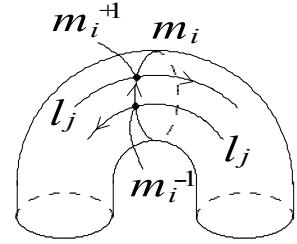

Fig. 5

A dual presentation from $(V; l, m)$ of $(U,V,F)$ is also defined in an analogous manner, and is denoted as

$$\pi_1(M^3) = \langle l_1,\cdots,l_n \mid \hat{m}_1 = 1,\cdots,\hat{m}_n = 1 \rangle \quad (1')$$

Group (1) is isomorphic to (1′) but the presentation is generally different from (1′) because meridians and longitudes are switched in $(U;m,l)$ and $(V;l,m)$. Therefore the forms of relators in (1) and (1′) are different generally.

Let a presentation of the fundamental group derived from U-A,U-B of $(U;m,l)$ be (UA),(UB), respectively.

$$\left\langle \begin{array}{c} m_i, m_j \\ m_k \\ (k \neq i,j) \end{array} \middle| \begin{array}{l} \underline{m_i m_j} w_{ij_k} = 1 \cdots (l_{ij_k})(k=1,\cdots,l) \\ m_i^{-1} w_{i_k} = 1 \cdots (l_{i_k})(k=1,\cdots,p) \\ m_j w_{j_k} = 1 \cdots (l_{j_k})(k=1,\cdots,q) \\ r_\alpha = 1 \text{(relations other than the above)} \end{array} \right\rangle \quad \text{(UA)}$$

$$\left\langle \begin{array}{c} m_i, m_j \\ m_k \\ (k \neq i,j) \end{array} \middle| \begin{array}{l} m_i w_{ij_k} = 1 \cdots (l_{ij_k})(k=1,\cdots,l) \\ m_j m_i^{-1} w_{i_k} = 1 \cdots (l_{i_k})(k=1,\cdots,p) \\ m_j w_{j_k} = 1 \cdots (l_{j_k})(k=1,\cdots,q) \\ r_\alpha = 1 \text{(relations other than the above)} \end{array} \right\rangle \quad \text{(UB)}$$

Note that the relations $(l_{ij_k})(k=1,\cdots,l)$ in (UA) and (UB), too, do not exist if the longitudes $(l_{ij_k})(k=1,\cdots,l)$ do not exist.

Operations; in (UA), replace the generator $m_i$, letters $m_i m_j$ in $(l_{ij_k})$ by a new letter $\tilde{m}_i$ (a new generator), we get a presentation that isomorphic to (UB).

Let a presentation of the fundamental group derived from V-A,V-B of $(V;l,m)$ be (VA),(VB), respectively.

$$\left\langle \begin{array}{l} l_{i_1},\cdots,l_{i_p} \\ l_{j_1},\cdots,l_{j_q} \\ l_{ij_1},\cdots,l_{ij_l} \\ l_k(k\neq i,j,ij) \end{array} \middle| \begin{array}{l} l_{i_1}^{-1}\cdots l_{i_p}^{-1}\underline{l_{ij_1}\cdots l_{ij_l}}=1\cdots(m_i) \\ l_{j_1}\cdots l_{j_q}\underline{l_{ij_1}\cdots l_{ij_l}}=1\cdots(m_j) \\ r_\alpha{'}=1\text{(relations other than the above)} \end{array} \right\rangle \quad \text{(VA)}$$

$$\left\langle \begin{array}{l} l_{i_1},\cdots,l_{i_p} \\ l_{j_1},\cdots,l_{j_q} \\ l_{ij_1},\cdots,l_{ij_l} \\ l_k(k\neq i,j,ij) \end{array} \middle| \begin{array}{l} l_{i_1}^{-1}\cdots l_{i_p}^{-1}l_{ij_1}\cdots l_{ij_l}=1\cdots(m_i) \\ l_{j_1}\cdots l_{j_q}l_{i_p}\cdots l_{i_1}=1\cdots(m_j) \\ r_\alpha{'}=1\text{(relations other than the above)} \end{array} \right\rangle \quad \text{(VB)}$$

Operation; in (VA), by substituting $l_{ij_1}\cdots l_{ij_l}=l_{i_p}\cdots l_{i_1}$ derived from $(m_i)$ for $l_{ij_1}\cdots l_{ij_l}$ in $(m_j)$, we get (VB).

In like manner, transformations of the fundamental groups corresponding to those of a handles sliding of $h_j$ about $h_i$ of $(U;m,l)$ and a band move of $m_i$ of $(V;l,m)$ are obtained. Therefore by gathering the Theorem 10 and considering the above, we have;

**Theorem 12.** *The transformation from U-A into U-B by a handle sliding of $h_i$ about $h_j$, or a band move of $m_j$ corresponds to the replacements letters of the fundamental group.*

*The dual transformation from V-A into V-B by a band move of $m_j$ corresponds to the substitution of the fundamental group.*

**Corollary 13.** *If the fundamental group derived from a Heegaard diagram of $M^3$ is transformed into the trivial group by a finite sequence of the replacement and substitution corresponding the handle sliding and band move in U-A,V-A, then $M^3$ is the 3-sphere.*

## References

[1] S. Horiguchi: Transformations of the fundamental groups corresponding to those of Heegaard diagrams by the band moves, Bulletin of Niigata Sangyo Univ., No.17, 155-164. 1997.6
   http://ci.nii.ac.jp/naid/110000484551

[2] H. Poincaré: Cinquième complément a 1'analysis situs, Rend. Circ. Mat. Palermo **18**(1904),

45-110.

[3] H. Seifert & W. Threlfa11: A Textbook of Topology, Translated in English by M.A. Goldman, Academic Press, lnc. 1980.

[4] J.S. Birman: Heegaard sp1ittlngs, diagrams and sewings for closed orientab1e 3-manifolds, Lecture notes for CBMS conference at Blacksburg, (1977).

[5] H. Zieschang: Über einfache Kurven auf Vollbrezeln, Abh. Math. Sem. Univ. Hamburg **25**(1962), 231-250.

[6] F. Waldhausen: Heegaard-Zerlegungen der 3-Sphäre, Topology **7**(1968) 195-203.

[7] J. Hempe1: 3-manifolds, Ann. of Math. Studies 86, Princeton Univ. Press. 1976.

Shunji HORIGUCHI  
Niigata Sangyo University  
4730, Karuigawa, Kashiwazaki,  
Niigata, 945-1393, Japan  
*E-mail*: shori@econ.nsu.ac.jp